\documentclass{article}
\usepackage {amsmath, amssymb, amsbsy, amsthm, graphicx, verbatim, eufrak, eucal}

\newtheorem{theorem}{Theorem}
\newtheorem{lemma}{Lemma}

\newtheorem{definition}{Definition}

\newtheorem{example}{Example}
\newtheorem{remark}{Remark}

\numberwithin{theorem}{section}
\numberwithin{definition}{section}
\numberwithin{lemma}{section}
\numberwithin{corollary}{section}
\numberwithin{equation}{section}
\numberwithin{proposition}{section}
\numberwithin{example}{section}
\numberwithin{remark}{section}
\numberwithin{figure}{section}

\def\vtl{\vskip 1mm}

\def\es{\varnothing}
\def\SB{\subseteq}

\def\aa{\alpha}
\def\bb{\beta}
\def\gg{\gamma}

\def\bm{\boldsymbol m}
\def\bn{\boldsymbol n}
\def\bp{\boldsymbol p}
\def\bq{\boldsymbol q}

\def\b0{\boldsymbol 0}

\def\Zee{\mathbb Z}

\def\imp{\Rightarrow}
\def\EQ{\Longleftrightarrow}
\def\eq{\Leftrightarrow}

\def\AAA{{\cal A}}
\def\FFF{{\cal F}}

\def\CCC{{\cal C}}

\def\HHH{{\cal H}}

\def\VVV{{\cal V}}

\def\ZZZ{{\cal Z}}

\def\PPP{{\cal P}}

\def\SSS{{\cal S}}
\def\TTT{{\cal T}}
\def\UUU{{\cal U}}

\def\BS{\bold S}

\def\BP{\bold P}

\def\POW{\mathfrak P}
\def\POWF{\mathfrak P_{\text{\sc f}}}

\def\begeq{\begin{equation}}
\def\edeq{\end{equation}}

\def\roster{\begin{enumerate}}
\def\endroster{\end{enumerate}}

\makeindex

\begin{document}

\title{Cubical Token Systems}
	\author{Sergei~Ovchinnikov\\ Mathematics Department\\San Francisco State University\\San Francisco, CA 94132\\sergei@sfsu.edu}

\date\empty
\maketitle

\begin{abstract}
\noindent
The paper deals with combinatorial and stochastic structures of cubical token systems. A cubical token system is an instance of a token system, which in turn is an instance of a transition system. It is shown that some basic results of combinatorial and stochastic parts of media theory hold almost in identical form for cubical token systems, although some underlying concepts are quite different. A representation theorem for a cubical token system is established asserting that the graph of such a system is cubical. 
\end{abstract}

\noindent
{\small{\sl Keywords}: Transition system, token system, cubical system, medium, Markov chain}

\section{Introduction} \label{S:intro}

Cubical token systems and media are particular instances of a general algebraic structure, called `token system', describing a mathematical, physical, or behavioral system as it evolves from one `state' to another. This structure is formalized as a pair $(\SSS,\TTT)$ consisting of a set $\SSS$ of states and a set $\TTT$ of tokens. Tokens are transformations of the set of states. Strings of tokens are `messages' of the token system. The concept of a medium was introduced in~\cite{jF97} as a token system specified by some constraining axioms, and developed further in~\cite{jF02,sO00}. For more recent advances in media theory the reader is referred to~\cite{sO06,sO06a} and the forthcoming monograph~\cite{dE07}.

In the field of computer science, tokens systems are special forms of `transition systems'~\cite{gW94}. However, we do not follow this lead in the paper. Instead, we propose a system of axioms specifying a class of token systems that we call `cubical (token) systems'. The name is justified by the result of Section~\ref{S:representation theorem} asserting that the graph of a cubical system is cubical.

We begin by introducing basic concepts of token systems in Section~\ref{S:defs} and axioms for cubical systems and media in Section~\ref{S:axioms}, where it is also shown that media form a subclass of cubical systems.

$G$-systems are token systems defined on connected families of sets. In Section~\ref{S:canonical example} we show that they are instances of cubical systems. $G$-systems are typical examples of cubical systems as it is demonstrated in Section~\ref{S:representation theorem}.

Structural properties of states and messages of a cubical system are established in Section~\ref{S:tokens and contents} in terms of their `contents'. These properties are crucial for the development of stochastic token theory presented in Section~\ref{S:stochastic}.

The main result of the algebraic part of cubical systems---the representation theorem---is established in Section~\ref{S:representation theorem} (Theorem~\ref{RT}).

In Section~\ref{S:examples}, we give some examples of cubical systems that could serve as potential applications.

\section{Token systems} \label{S:defs}

Let $\SSS$ be a set of {\em states}. A {\em token} is a transformation $\tau:S\mapsto S\tau$. By definition, the identity function $\tau_0$  on $\SSS$ is not a token. Let $\TTT $ be a set of tokens. The pair $(\SSS,\TTT )$ is called a {\em token system}. To avoid trivialities, we assume that  $|\SSS|\geq 2$ and $\TTT\neq \es$.

Let $V$ and $S$ be two states of a token system $(\SSS,\TTT)$. Then $V$ is {\em adjacent} to $S$ if $S\neq V$ and $S\tau = V$ for some token $\tau\in\TTT$. A token $\tilde \tau$ is a {\em
reverse} of a token $\tau$ if for all distinct
$S,V\in\SSS$, we have 
$$
S\tau=V \quad\EQ\quad  V\tilde \tau=S.
$$
Two distinct states $S$ and $V$ are {\em adjacent} if $S$ is adjacent to $V$ and $V$ is adjacent to $S$.

\begin{remark} \label{token remarks}
{\rm Suppose that tokens $\mu$ and $\nu$ are reverses of a token $\tau$. Then, for $S\neq V$,
$$
V\mu=S\quad\eq\quad S\tau=V\quad\eq\quad V\nu=S.
$$
It follows that $\mu=\nu$. Therefore, if a reverse of a token exists, then it is unique. It is also clear that the reverse of a reverse is the token itself, $\tilde{\tilde\tau}=\tau$, provided that $\tilde\tau$ exists.

{\begin{figure}[h!]
\centerline{\includegraphics{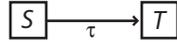}}
\caption{Token system $(\SSS,\TTT)$ with $\SSS=\{S,T\}$ and $\TTT=\{\tau\}$. The state $T$ is adjacent to state $S$ but these two states are not adjacent in $(\SSS,\TTT)$.} \label{trivial-TS} 
\end{figure}
}
In general, a token of a token system $(\SSS,\TTT)$ does not necessarily have a reverse in $(\SSS,\TTT)$. For instance, the token $\tau$ of the token system shown in Figure~\ref{trivial-TS} does not have a reverse in $\TTT$. It is also possible for a token to be a reverse of itself (see Example [C1] in Figure~\ref{4-axioms}).
}
\end{remark}

A {\em message} of a token system $(\SSS,\TTT)$ is a string of elements of the set $\TTT$. We write these strings in the form $\bm=\tau_1\tau_2\ldots\tau_n$. If a token $\tau$ occurs in the string $\tau_1\tau_2\ldots\tau_n$, we say that the message $\bm=\tau_1\tau_2\ldots\tau_n$ {\em contains} $\tau$. 

A message $\bm=\tau_1\tau_2\ldots\tau_n$ defines a transformation
$$
S\mapsto S\bm=((\ldots((S\tau_1)\tau_2)\ldots)\tau_n)
$$
of the set of states $\SSS$. By definition, the empty message defines the identity transformation $\tau_0$ of $\SSS$. If $V=S\bm$ for some message $\bm$ and states $S,V\in\SSS$, then we say that $\bm$ {\em produces} $V$ from $S$ or, equivalently, that $\bm$ transforms $S$ into $V$. More generally, if $\bm=\tau_1\ldots\tau_n$, then we say that $\bm$ {\em produces} a sequence of states $(S_i)$, where $S_i=S\tau_0\tau_1\ldots\tau_i$ for $0\le i\le n$.

If $\bm$ and $\bn$ are two messages, then $\bm\bn$ stands for the concatenation of the strings $\bm$ and $\bn$. We denote by $\widetilde\bm=\tilde\tau_n\ldots\tilde\tau_1$ the {\em reverse} of the message $\bm=\tau_1\ldots\tau_n$, provided that the tokens in $\widetilde\bm$ exist.

A message $\bm=\tau_1\ldots\tau_n$ is {\em vacuous} if the set of indices $\{1,\ldots,n\}$ can be partitioned into pairs $\i,j$ with $i\neq j$, such that $\tau_i$ and $\tau_j$ are mutual reverses.

A~message $\bm$ is {\em effective} (respectively {\em ineffective})  for a state $S$ if $S\bm\neq S$ (respectively $S\bm=S$) for the corresponding transformation $\bm$. A~message $\bm=\tau_1\ldots\tau_n$ is {\em stepwise effective} for $S$ if  $S_k\neq S_{k-1}$, $1\leq k\leq n$, in the sequence of states produced by $\bm$ from $S$. A message is {\em closed} for a state $S$ if it is stepwise effective and ineffective for $S$. When it is clear from the context which state is under consideration, we may drop a reference to that state.

Two token systems $(\SSS,\TTT)$ and $(\SSS',\TTT')$ are said to be {\em isomorphic} if there is a pair $(\aa,\bb)$ of bijections $\aa:\SSS\rightarrow\SSS'$ and $\bb:\TTT\rightarrow\TTT'$ such that
$$
S\tau=T\quad\eq\quad \aa(S)\bb(\tau)=\aa(T)
$$
for all $S,T\in\SSS$ and $\tau\in\TTT$.

\section{Axioms for cubical systems} \label{S:axioms}

\begin{definition}
{\rm A token system $(\SSS,\TTT)$ is called a {\em cubical (token) system} if the following axioms are satisfied:
\roster
\item[]
	\roster
	\item[{[C1]}] Every token $\tau\in\TTT$ has a reverse $\tilde\tau\in\TTT$ and $\tilde\tau\neq\tau$.
	\item[{[C2]}] For any two distinct states $S$ and $T$ there is a stepwise effective message producing $T$ from $S$.
	\item[{[C3]}] A message which is stepwise effective for some state is closed for that state if and only if it is vacuous.
	\item[{[C4]}] If $\bm=\tau_1\ldots\tau_n$ is a stepwise effective message for some state, then occurrences of a token and its reverse alternate in $\bm$. More specifically, if $\tau_i=\tau_j=\tau$ for $i<j$ and some $\tau\in\TTT$, then $\tau_k=\tilde\tau$ for some $i<k<j$.
	\endroster
\endroster
}
\end{definition}

\begin{theorem}
Axioms {\rm[C1]--[C4]} are independent.
\end{theorem}

\begin{proof}
Each diagram in Figure~\ref{4-axioms} shows a token system satisfying exactly three of the four axioms defining  cubical systems. Each drawing is labeled by the failing axiom. We omit the proofs.
\end{proof}

{\begin{figure}[h!]
\centerline{\includegraphics{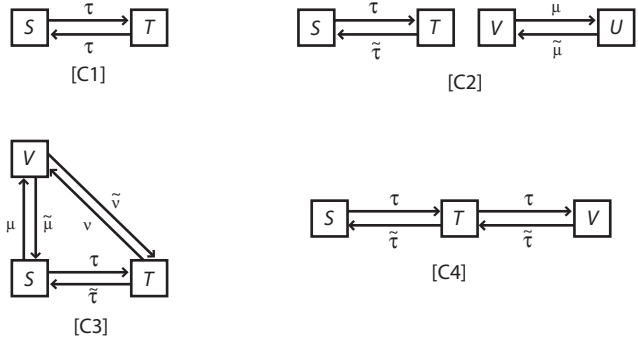}}
\caption{Digraphs of four token systems. Loops are omitted. Each system is labeled by the unique failing axiom.} \label{4-axioms} 
\end{figure}
}

We need the concept of a `concise message' for the definition of a medium.

\begin{definition}
{\rm A message $\bm$ is said to be {\em concise} for a state $S$ if: (i) $\bm$ is stepwise effective for $S$, (ii) no token occurs twice in $\bm$, and $\bm$ does not contain a token and its reverse.
}
\end{definition}

\begin{definition}
{\rm A token system $(\SSS,\TTT)$ is called a {\em medium} (on $\SSS$) if the following axioms are satisfied.
\roster
\item[]
\roster
\item[{[Ma]}] For any two distinct states $S$ and $V$ in $\SSS$ there is a
concise message transforming $S$ into $V$.
\item[{[Mb]}] A message which is closed for some state is vacuous.
\endroster
\endroster
}
\end{definition}

\begin{theorem}
A medium is a cubical system.
\end{theorem}

\begin{proof}
Let $(\SSS,\TTT)$ be a medium. Axiom [C1] follows from Lemma~5.1 in~\cite{sO06a}, Axiom [C2] is an immediate consequence of Axiom [Ma], and Axiom [C3] follows from [Mb] and Lemma 5.4 in~\cite{sO06a}.

It remains to verify that Axiom [C4] holds for $(\SSS,\TTT)$. Let $\bm=\tau_1\ldots\tau_n$ be a stepwise effective message for a state $S$ and $(S_i)$ be a sequence of states produced by $\bm$. Suppose that $\tau_i$ and $\tau_j$ ($i<j$) are two consecutive occurences of a token $\tau$ in $\bm$ such that there is no occurrence of $\tilde\tau$ between $\tau_i$ and $\tau_j$. By [Ma], there is a concise message $\bn$ producing $S_{i-1}$ from $S_j$. By [Mb], we must have two occurences of $\tilde\tau$ in the concise message $\bn$, a contradiction. Thus [C4] holds for a medium.
\end{proof}

{\begin{figure}[h!]
\centerline{\includegraphics{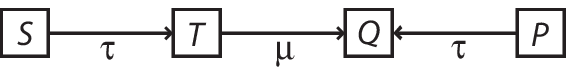}}
\caption{Token system $(\SSS,\TTT)$ with $\SSS=\{S,T,P,Q\}$ and $\TTT=\{\tau,\tilde\tau,\mu,\tilde\mu\}$.} \label{QM example} 
\end{figure}
}

\begin{example}
{\rm Let $(\SSS,\TTT)$ be a token system displayed in Figure~\ref{QM example}. There is no concise message producing $P$ from $S$, so this token system is not a medium. It is easy to verify that this system is a cubical system.
}
\end{example}

\section{A `canonical' example of a cubical system} \label{S:canonical example}

A `canonical' example of a medium is the representing medium of a well-graded family of sets~\cite{jF02,sO06a,sO00}. For cubical systems, similar examples are given by $G$-systems.

\begin{definition}
{\rm A {\em cube $\HHH(X)$ on a set $X$} is a graph that has the set of all finite subsets of $X$ as the set of vertices; $\{S,T\}$ is an edge of $\HHH(X)$ if $|S\bigtriangleup T|=1$. A graph is said to be {\em cubical} if it is embeddable into some cube $\HHH(X)$. A {\em partial cube} is a graph that is isometrically embeddable into some cube $\HHH(X)$. 
}
\end{definition}

\begin{definition}
{\rm Let $G=(\FFF,\VVV)$ be a connected subgraph of the cube $\HHH(X)$ on a set $X$ with $|X|\ge 2$. A {\em $G$-system on $\FFF$} is a pair $(\FFF,\TTT_G)$ where $\TTT_G$ is a family of transformations defined by
\begeq \label{gamma}
\gg_x : S\mapsto S\gg_x = \begin{cases}
	S\cup\{x\}, &\text{if $\{S,S\cup\{x\}\}\in\VVV$,}\\
	S, &\text{otherwise,}
\end{cases}
\edeq
\begeq \label{tilde-gamma}
\tilde\gg_x : S\mapsto S\tilde\gg_x = \begin{cases}
	S\setminus\{x\}, &\text{if $\{S,S\setminus\{x\}\}\in\VVV$,}\\
	S, &\text{otherwise,}
\end{cases}
\edeq
for $x\in\cup\,\FFF\setminus\cap\,\FFF$.
}
\end{definition}

\begin{example}
{\rm Let $X=\{x,y\}$ and $G=(\FFF,\VVV)$, where, $\FFF=\POW(X)$ and 
$$\VVV=\{\{\es,\{x\}\},\{\{x\},\{x,y\}\},\{\{y\},\{x,y\}\}\}.$$
This $G$-system on $\FFF$ is isomorphic to the token system displayed in Figure~\ref{QM example} under isomorphism $(\aa,\bb)$ defined by
\begin{gather*}
\aa(S)=\es,\quad\aa(T)=\{x\},\quad\aa(Q)=\{x,y\},\quad\aa(P)=\{y\},\\
\bb(\tau)=\gg_x,\quad\bb(\tilde\tau)=\tilde\gg_x,\quad\bb(\mu)=\gg_y,\quad\bb(\tilde\mu)=\tilde\gg_y.
\end{gather*}
Note that $\es\gg_y=\es\neq\{y\}$.
}
\end{example}

\begin{theorem}
A $G$-system on $\FFF$ is a token system and, for any $x\in\cup\,\FFF\setminus\cap\,\FFF$, the tokens $\gg_x$ and $\tilde\gg_x$ are mutual reverses.
\end{theorem}

\begin{proof}
We show first that the functions defined by~(\ref{gamma}) and~(\ref{tilde-gamma}) are tokens. Clearly, for any $x\in\cup\,\FFF\setminus\cap\,\FFF$ there are two sets $S,T\in\FFF$ such that $x\notin S$ and $x\in T$. Since $G$ is a connected graph, there is a sequence $S_0=S,S_1,\ldots,S_n=T$ of sets in $\FFF$ such that $\{S_i,S_{i+1}\}\in\VVV$ for $0\le i<n$. In particular, $|S_i\bigtriangleup S_{i+1}|=1$. Since $x\notin S$ and $x\in T$, there is $j$ such that $x\notin S_j$ and $x\in S_{j+1}$, so, by~(\ref{gamma}),
$$
S_{j+1}=S_j\cup\{x\}=S_i\gg_x.
$$
It follows that $\gg_x\neq\tau_0$ for any $x\in\cup\,\FFF\setminus\cap\,\FFF$. The case of functions defined by~(\ref{tilde-gamma}) is treated similarly. It is clear that the tokens $\gg_x$ and $\tilde\gg_x$ are mutual reverses.
\end{proof}

Let $S_0=S,S_1,\ldots,S_n=T$ be a walk in $G$. For an edge $\{S_{i-1},S_i\}$, we denote $\{x_i\}=S_{i-1}\bigtriangleup S_i$, $\tau_i=\gg_{x_i}$ if $S_i=S_{i-1}\cup\{x_i\}$, and $\tau_i=\tilde\gg_{x_i}$, otherwise. Then $\bm=\tau_1\ldots\tau_n$ is a stepwise effective message for $S$ of the $G$-system $(\FFF,\TTT_G)$. Conversely, a stepwise effective message $\bm=\tau_1\ldots\tau_n$ of $(\FFF,\TTT_G)$ producing a state $T$ from a state $S$ defines a walk $W_{\bm}$ in $G$ with vertices $S_i=S\tau_0\tau_1\ldots\tau_n$. Thus there is a one-to-one correspondence between the stepwise effective messages of a $G$-system and the walks in $G$.

\begin{theorem}
A $G$-system on $\FFF$ is a cubical system on the set of states $\FFF$.
\end{theorem}

\begin{proof}
Let $(\FFF,\TTT_G)$ be a $G$-system on $\FFF$. Axiom [C1] holds trivially and [C2] holds because $G$ is a connected graph.

Let $\bm=\tau_1\ldots\tau_n$ be a stepwise effective message for a state $S$. Suppose that there are two consecutive occurences of $\gg_x$ in $\bm$, say, $\tau_i=\gg_x$ and $\tau_j=\gg_x$ with $i<j$, such that there is no occurrence of $\tilde\gg_x$ between $\tau_i$ and $\tau_j$. Then $x\in S_i=S\tau_0\tau_1\ldots\tau_i$ which implies $x\in S_{j-1}$, since $\tilde\gg_x$ does not occur between $\tau_i$ and $\tau_j$. It follows that $\gg_x=\tau_j$ is not effective for the state $S_{j-1}$, a contradiction. Thus occurences of a token and its reverse must alternate in $\bm$, so [C4] holds for $(\FFF,\TTT_G)$. A minor modification of this argument shows that [C3] also holds for $(\FFF,\TTT_G)$.
\end{proof}

\section{Tokens and contents} \label{S:tokens and contents}

Tokens of a cubical system share many properties with tokens of a medium (cf. Lemmas 5.1 and 5.2 in~\cite{sO06a}).

\begin{lemma} \label{tokens}
The following statements hold for a cubical system $(\SSS,\TTT)$:
\roster
	\item[{\rm(i)}] $\tilde{\tilde\tau}=\tau$ for any $\tau\in\TTT$.
	\item[{\rm(ii)}] For any two adjacent states $S$ and $T$ there is a unique token producing $T$ from $S$.
	\item[{\rm(iii)}] If $S$, $T$, and $P$ are three distinct states such that $S\tau=T$ and $T\mu=P$, for some tokens $\tau$ and $\mu$, then $\mu\neq\tau$ and $\mu\neq\tilde\tau$.
	\item[{\rm(iv)}] No token can be a one-to-one function.
\endroster
\end{lemma}

\begin{proof}
(i) By [C1], $\tilde\tau$ exists, so $\tilde{\tilde\tau}=\tau$ (cf. Remark~\ref{token remarks}).
\vtl
(ii) Suppose that $S\tau=S\mu=T$. By [C1] and [C3], the message $\tau\tilde\mu$ is well-defined and vacuous, so $\tau=\tilde{\tilde\mu}=\mu$.
\vtl
(iii) Since $\tau\mu$ is a stepwise effective message for $S$, we have $\mu\not=\tau$, by [C4]. If $\mu=\tilde\tau$, then $S=T\tilde\tau=P$, a contradiction, since $S\neq P$ and $\tilde\tau$ is a function.
\vtl
(iv) Since $\tau$ is not the identity transformation, there are states $S$ and $T$ such that $S\tau=T$. By (iii), $S\tau=T=T\tau$, so $\tau$ is not a one-to-one function.
\end{proof}

\begin{remark}
{\rm Property (ii) of Lemma~\ref{tokens} is a very strong property of tokens of a cubical system. It asserts that two tokens $\tau$ and $\mu$ transforming some state $S$ into a different state $T$ are equal transformations, that is, $V\tau=V\mu$ for all $V\in\SSS$.
}
\end{remark}

Let $\tau$ be a token of a medium. We define
$$
\UUU_\tau=\{S\in\SSS\mid S\tau\neq S\}.
$$
Note that $\UUU_\tau\neq\es$, since $\tau$ is a token.

\begin{lemma} \label{token-structure}
For any given $\tau\in\TTT$ we have
\roster
	\item[{\rm(i)}] $(\UUU_\tau)\tau=\UUU_{\tilde\tau}$.
	\item[{\rm(ii)}] $\UUU_\tau\cap\UUU_{\tilde\tau}=\es$.
	\item[{\rm(iii)}] The restriction $\tau|_{_{\UUU_\tau}}$ is a bijection from $\UUU_\tau$ onto $\UUU_{\tilde\tau}$ with $\tau|^{-1}_{_{\UUU_\tau}}=\tilde\tau|_{_{\UUU_{\tilde\tau}}}$.
\endroster
\end{lemma}

\begin{proof}
(i) We have
$$
T\in(\UUU_\tau)\tau\;\eq\; S\tau=T\;(S\neq T)\;\eq\; T\tilde\tau=S\;(S\neq T)\;\eq\; T\in\UUU_{\tilde\tau}.
$$

(ii) If $S\in\UUU_\tau\cap\UUU_{\tilde\tau}$, then there exist $T\neq S$ such that $S\tau=T$, and $V\neq S$ such that $S\tilde\tau=V$, so $V\tau=S$. By Lemma~\ref{tokens}(ii) and Axiom [C1], $V\neq T$ contradicting Lemma~\ref{tokens}(iii). It follows that $\UUU_\tau\cap\UUU_{\tilde\tau}=\es$.

(iii) follows immediately from (i) and (ii).
\end{proof}

\begin{definition}
{\rm Let $(\SSS,\TTT)$ be a cubical system. For any token $\tau$ and any message $\bm$, we define $\#(\tau,\bm)$ as the number of occurrences of $\tau$ in $\bm$.
For any message $\bm$, the {\em content} of $\bm$ is the set $\CCC(\bm)$ defined by
$$
\CCC(\bm)=\{\tau\in\TTT\mid \#(\tau,\bm)>\#(\tilde\tau,\bm)\}.
$$ 
For any state $S$, the {\em content} $\widehat{S}$ of $S$ is the union $\cup_{\bm}\CCC(\bm)$ taken over the set of all stepwise effective messages producing the state $S$.
}
\end{definition}

The two concepts of `content' are different from their counterparts in media theory. For instance, the content of a vacuous message of a cubical system is empty, whereas it is not empty in media theory. However, the main results of media theory concerning these concepts are valid for cubical systems. We establish these results in a series of theorems in the rest of this section. Note that the results of Theorems~\ref{hatV-hatS} and~\ref{hatS=hatV} are especially useful in stochastic part of cubical systems theory (Section~\ref{S:stochastic}). 

The following properties of the functions $\#$ and $\CCC$ are immediate and will be used implicitly in the paper:
\begin{gather*}
\#(\tilde\tau,\widetilde\bm)=\#(\tau,\bm),\quad \#(\tilde\tau,\bm)=\#(\tau,\widetilde\bm),\\
\tau\in\CCC(\bm)\;\eq\;\tilde\tau\in\CCC(\widetilde\bm),\quad\tau\in\CCC(\bm)\;\imp\;\tilde\tau\notin\CCC(\bm).
\end{gather*}

\begin{lemma}
If $\bm$ is a stepwise effective message for some state, then
\begeq \label{n=n+1}
\tau\in\CCC(\bm)\quad\eq\quad \#(\tau,\bm)=\#(\tilde\tau,\bm)+1.
\edeq
Therefore, for any $\tau\in\TTT$,
\begeq \label{n-n}
\#(\tau,\bm)-\#(\tilde\tau,\bm)\in\{-1,0,1\}.
\edeq
\end{lemma}

\begin{proof}
By Axiom [C4] the occurences of $\tau$ and $\tilde\tau$ in $\bm$ alternate. Therefore,
$$
\tau\in\CCC(\bm)\quad\imp\quad\#(\tau,\bm)>\#(\tilde\tau,\bm)\quad\imp\quad\#(\tau,\bm)=\#(\tilde\tau,\bm)+1.
$$
The converse implication in~(\ref{n=n+1}) is trivial. It is clear, that~(\ref{n-n}) follows from~(\ref{n=n+1}).
\end{proof}

\begin{lemma} \label{no token+reverse}
The content of a state cannot contain both a token and its reverse.
\end{lemma}

\begin{proof}
Suppose that $\tau,\tilde\tau\in\widehat S$ for some token $\tau$ and some state $S$. Then there are two stepwise effective messages $\bm$ and $\bn$ both producing $S$ and such that $\tau\in\CCC(\bm)$ and $\tilde\tau\in\CCC(\bn)$. Therefore,
$$
\#(\tau,\bm)=\#(\tilde\tau,\bm)+1\quad\text{and}\quad\#(\tilde\tau,\bn)=\#(\tau,\bn)+1.
$$
It follows that
\begin{align*}
\#(\tau,\bm\widetilde\bn)&=\#(\tau,\bm)+\#(\tau,\widetilde\bn)=\#(\tilde\tau,\bm)+1+\#(\tilde\tau,\widetilde\bn)+1 \\
&=\#(\tilde\tau,\bm\widetilde\bn)+2,
\end{align*}
which contradicts Axiom [C4] since $\bm\widetilde\bn$ is a stepwise effective message for some state. Therefore $\widehat S$ cannot contain both $\tau$ and $\tilde\tau$.
\end{proof}

\begin{theorem}
For any token $\tau$ and any state $S$ of a cubical system, we have either $\tau\in\widehat S$ or $\tilde\tau\in\widehat S$ (but not both).
\end{theorem}

\begin{proof}
Since $\tau$ is a token, there are distinct states $V$ and $W$ such that $V\tau=W$. By Axiom [C2], there are stepwise effective messages $\bm$ and $\bn$ such that $S\bm=V$ and $W\bn=S$. (If $S$ equals either $V$ or $W$, the corresponding message is empty.) By Axiom [C3], the message $\tau\bn\widetilde\bm$ is vacuous. Therefore,
$$
\#(\tau,\tau\bn\widetilde\bm)=\#(\tilde\tau,\tau\bn\widetilde\bm).
$$
We have
$$
\#(\tau,\tau\bn\widetilde\bm)=1+\#(\tau,\bn)+\#(\tilde\tau,\bm)
$$
and
$$
\#(\tilde\tau,\tau\bn\widetilde\bm)=\#(\tilde\tau,\bn)+\#(\tau,\bm).
$$
From the last three displayed equations we obtain
$$
[\#(\tau,\bm)-\#(\tilde\tau,\bm)]+[\#(\tilde\tau,\bn)-\#(\tau,\bn)]=1.
$$
By~(\ref{n-n}), we must have either $\#(\tau,\bm)-\#(\tilde\tau,\bm)=1$ or $\#(\tilde\tau,\bn)-\#(\tau,\bn)=1$ but not both. It follows that either $\tau\in\CCC(\bm)$ or $\tilde\tau\in\CCC(\bn)$. By Lemma~\ref{no token+reverse}, either $\tau\in\widehat S$ or $\tilde\tau\in\widehat S$.
\end{proof}

\begin{theorem} \label{hatV-hatS}
If $S$ and $V$ are two distinct states, with $S\bm=V$ for some stepwise effective message $\bm$, then $\widehat V\setminus\widehat S=\CCC(\bm)$. Therefore,
$$
\widehat S\bigtriangleup\widehat V=\CCC(\bm)+\CCC(\widetilde\bm),
$$
where $+$ stand for the disjoint union of two sets. In particular,
$$
\widehat S\bigtriangleup\widehat V=\{\tau,\tilde\tau\},
$$
if $S\tau=V$.
\end{theorem}

\begin{proof}
Let $\tau$ be a token in $\CCC(\bm)$. Then $\tau\in\widehat V$ and $\tilde\tau\in\CCC(\widetilde\bm)$ implying that $\tilde\tau\in\widehat S$. By Lemma~\ref{no token+reverse}, $\tau\notin\widehat S$. It follows that $\tau\in \widehat V\setminus\widehat S$. Thus $\CCC(\bm)\SB\widehat V\setminus\widehat S$.

Suppose now that $\tau\in \widehat V\setminus\widehat S$, so $\tau\in\widehat V$ and $\tau\notin\widehat S$. There is a stepwise effective message $\bn$ producing $V$ and such that $\tau\in\CCC(\bn)$. By~(\ref{n=n+1}),
\begeq \label{n-n=1}
\#(\tau,\bn)-\#(\tilde\tau,\bn)=1.
\edeq
Since $\tau\notin\widehat S$, we have $\tau\notin\CCC(\widetilde\bm)$ which implies, by~(\ref{n=n+1}) and~(\ref{n-n}),
\begeq \label{n-n<1}
\#(\tau,\widetilde\bm)-\#(\tilde\tau,\widetilde\bm)\in\{-1,0\}.
\edeq
The message $\bn\widetilde\bm$ is stepwise effective and produces the state $S$. We have $\tau\notin\CCC(\bn\widetilde\bm)$, since $\tau\notin\widehat S$. Therefore, by~(\ref{n=n+1}),~(\ref{n-n}) and~(\ref{n-n=1}),
\begin{align*}
1&>\#(\tau,\bn\widetilde\bm)-\#(\tilde\tau,\bn\widetilde\bm)=[\#(\tau,\bn)+\#(\tau,\widetilde\bm)]-[\#(\tilde\tau,\bn)+\#(\tilde\tau(\widetilde\bm)]\\
&=[\#(\tau(\bn)-\#(\tilde\tau,\bn)]+[\#(\tau,\widetilde\bm)-\#(\tilde\tau,\widetilde\bm)]=1+[\#(\tau,\widetilde\bm)-\#(\tilde\tau,\widetilde\bm)].
\end{align*}
By~(\ref{n-n<1}), $\#(\tilde\tau,\widetilde\bm)-\#(\tau,\widetilde\bm)=1$, or, equivalently, $\#(\tau,\bm)-\#(\tilde\tau,\bm)=1$. By~(\ref{n=n+1}), $\tau\in\CCC(\bm)$. Hence, $\widehat V\setminus\widehat S\SB\CCC(\bm)$. The result follows.
\end{proof}

\begin{lemma} \label{closed msg}
A stepwise effective message $\bm$ is closed if and only if $\CCC(\bm)=\es$.
\end{lemma}

\begin{proof}
A closed stepwise effective message $\bm$ is vacuous by Axiom [C3]. By~(\ref{n=n+1}), $\CCC(\bm)=\es$.

Conversely, if $\CCC(\bm)=\es$ for some stepwise effective message $\bm$, then, by Axiom [C4] and~(\ref{n=n+1}), $\bm$ is vacuous. By Axiom [C3], $\bm$ is closed.
\end{proof}

\begin{theorem} \label{hatS=hatV}
For any two states $S$ and $V$ we have
$$
\widehat S=\widehat V\quad\eq\quad S=T.
$$
\end{theorem}

\begin{proof}
Suppose that $\widehat S=\widehat V$ and let $\bm$ be a stepwise effective message producing $V$ from $S$. By Theorem~\ref{hatV-hatS}, $\CCC(\bm)=\es$. By Lemma~\ref{closed msg}, $\bm$ is a closed message. Thus, $S=V$. The converse implication is trivial.
\end{proof}

\begin{theorem}
Let $\bm$ and $\bn$ be two stepwise effective messages transforming some state $S$. Then
$$
S\bm=S\bn\quad\eq\quad\CCC(\bm)=\CCC(\bn).
$$
\end{theorem}

\begin{proof}
Suppose that $S\bm=S\bn=V$. By Theorem~\ref{hatV-hatS}, 
$$
\CCC(\bm)=\widehat V\setminus\widehat S=\CCC(\bn).
$$
Conversely, suppose that $\CCC(\bm)=\CCC(\bn)$ and let $V=S\bm$ and $W=S\bn$. By Theorem~\ref{hatV-hatS},
$$
\widehat V\bigtriangleup\widehat S=\CCC(\bm)+\CCC(\widetilde\bm)=\CCC(\bn)+\CCC(\widetilde\bn)=\widehat W\bigtriangleup\widehat S,
$$
implying $\widehat V=\widehat W$. By Theorem~\ref{hatS=hatV}, $V=W$.
\end{proof}

We conclude this section by comparing two concepts of contents with their counterparts in media theory.

\begin{theorem}
{\rm(i)} If $\bm=\tau_1\ldots\tau_n$ is a concise message of a medium, then
$$
\CCC(\bm)=\{\tau_1,\ldots,\tau_n\}.
$$

{\rm(ii)} For any state $S$ of a medium, its content $\widehat S$ is the set of all tokens each of which is contained in at least one concise message producing $S$.
\end{theorem}

\begin{proof}
As the first statement of the proposition is trivial, we proceed with a proof of (ii). If $\tau\in\CCC(\bm)$ for some concise message $\bm$ producing $S$, then, clearly, $\tau\in\widehat S$. Conversely, let $\tau\in\widehat S$ and $\bm$ be a stepwise effective message producing $S$ from some state $V$ and such that $\tau\in\CCC(\bm)$. By Axiom [Ma], there is a concise message $\bn$ producing $S$ from $T$. By~(\ref{n=n+1}), $\#(\tau,\bm)=\#(\tilde\tau,\bm)+1$. Therefore, by Axiom [Mb], $\tau\in\CCC(\bn)$. The result follows.
\end{proof}

\section{A representation theorem for cubical systems} \label{S:representation theorem}

\begin{definition}
{\rm The {\em graph} $G$ of a cubical system $(\SSS,\TTT)$ has $\SSS$ as the set of its vertices; two vertices are adjacent in $G$ if the corresponding states are adjacent in $(\SSS,\TTT)$.
}
\end{definition}

\begin{theorem} \label{RT}
Let $(\SSS,\TTT)$ be a cubical system. There exists a connected subgraph $G=(\FFF,\VVV)$ of some cube $\HHH(X)$ such that $(\SSS,\TTT)$ is isomorphic to the $G$-system $(\FFF,\TTT_G)$ on the family $\FFF$.
\end{theorem}

\begin{proof}
By Axiom [C2], the graph $G$ of the cubical system $(\SSS,\TTT)$ is connected. Let $J=\{\{\tau,\tilde\tau\}\}_{\tau\in\TTT}$. Elements of $J$ are called {\em labels}. By Lemma~\ref{tokens}(ii), a unique label is assigned to each edge of $G$.

We begin by constructing the family $\FFF$.

Let $S_0$ be a fixed state of the cubical system $(\SSS,\TTT)$. By [C2], for any state $T\neq S_0$, there is a stepwise effective message $\bm$ such that $S_0\bm=T$. We denote $W_{\bm}$ the walk in $G$ produced by the message $\bm$ and define a set $J_T$ by
$$
J_T=\{j\in J\mid \text{the number of occurrences of $j$ in $W_{\bm}$ is odd}\}.
$$
By definition, $J_{S_0}=\es$.

We need to show that the sets $J_T$ are well-defined. Suppose that $\bn$ is another stepwise effective message producing $T$ from $S_0$. By Axiom [C3], the number of occurrences of $j\in J$ in the closed walk $W_{\bm}W_{\widetilde\bn}$ is even. Hence, the number of occurrences of $j$ in $W_{\bm}$ is odd if and only if the number of its occurrences in $W_{\bn}$ is odd. Thus the set $J_T$ is well-defined and the assignment $T\mapsto J_T$ defines a mapping $\aa:\SSS\rightarrow\POWF(J)$, where $\POWF(J)$ stands for the family of finite subsets of $J$.

Let us prove that $\aa$ is a one-to-one mapping. Let $J_S=J_T$ for some states $S$ and $T$. By Axiom [C2] there are stepwise effective messages $\bm$, $\bn$, and $\bp$ such that $S_0\bm=S$, $S\bn=T$, and $T\bp=S_0$, so $W_{\bm}W_{\bn}W_{\bp}$ is a closed walk in $G$. By Axiom [C3], any label $j\in J$ occurs an even number of times in this walk. If $j\in J_S=J_T$, then $j$ occurs an odd number of times in each walks $W_{\bm}$ and $W_{\bp}$. Hence, $j$ occurs an even number of times in $W_{\bn}$. If $j\notin J_S=J_T$, then $j$ occurs an even number of times in each walks $W_{\bm}$ and $W_{\bp}$. Hence, $j$ occurs an even number of times in $W_{\bn}$. Thus any label occurs an even number of times in $W_{\bn}$. By Axiom [C4], the message $\bn$ is vacuous, and, by Axiom [C3], $S=T$. Hence, $\aa$ is a one-to-one mapping.

We show now that $\aa$ is an embedding of $G$ into the cube $\HHH(J)$. The sets $J_S$'s are vertices of the cube $\HHH(J)$. Let $P$ and $Q$ be two adjacent states of the cubical system $(\SSS,\TTT)$, so $P\tau=Q$ for some $\tau\in\TTT$, and let $j=\{\tau,\tilde\tau\}$ be the label of the edge $\{P,Q\}$ in the graph $G$. By Axiom [C2], there are stepwise effective messages $\bp$ and $\bq$ producing states $P$ and $Q$, respectively, from $S_0$. By Axiom [C3], $j$ occurs an even number of times in the closed walk $W_{\bp}W_{\tau}W_{\,\widetilde\bq}$. It follows that the label $j$ occurs an odd number of times either in $W_{\bp}$ or in $W_{\bq}$, so $j\in J_P\bigtriangleup J_Q$. Any other label $k$ occurs an even number of times in the walk $W_{\,\widetilde\bp}W_{\bq}$, so $k\notin J_P\bigtriangleup J_Q$. Thus, $J_P\bigtriangleup J_Q=\{j\}$, so $\{\aa(P),\aa(Q)\}$ is an edge of $\HHH(J)$. It follows that $\aa$ defines an embedding of the graph $G$ into the cube $\HHH(J)$. In the rest of the proof we identify $\aa(G)$ with $G$.

Let $\FFF=\{J_S\}_{S\in\SSS}$ and $(\FFF,\TTT_G)$ be the corresponding $G$-system on $\FFF$. 
(Clearly, $\cap\,\FFF=\es$ and $\cup\,\FFF=J$.) We prove that the cubical system $(\SSS,\TTT)$ is isomorphic to $(\FFF,\TTT_G)$.

Let $P$ and $Q$ be two adjacent states of the cubical system $(\SSS,\TTT)$. Since $\aa(P)=J_P$ and $\aa(Q)=J_Q$ are adjacent in the graph $G$, we have $J_P\bigtriangleup J_Q=\{j\}$ for some $j\in J$, so we may assume that $J_Q=J_P+\{j\}$. Since $P$ and $Q$ are adjacent states, we have $P\tau=Q$ for some token $\tau$. Note that $j=\{\tau,\tilde\tau\}$. We define $\bb(\tau)=\gg_j$, $\bb(\tilde\tau)=\tilde\gg_j$ and show that these assignments do not depend on a particular choice of $P$ and $Q$ with $P\tau=Q$. Let $S$ and $T$ be another pair of adjacent states such that $S\tau=T$, and let $\bm$ and $\bn$ be stepwise effective messages producing $Q$ from $S_0$ and $T$ from $Q$, respectively. By Axiom [C4], there is an even number of occurrences of the label $j$ in the walk $W_\tau W_{\bn} W_{\tilde\tau}$ connecting $P$ with $S$, so there is an even number of occurrences of $j$ in $W_{\bn}$. Since $j\in J_Q$, there is an odd number of occurrences of $j$ in $W_{\bm}$. Therefore, there is an odd number of occurrences of $j$ in the walk $W_{\bm}W_{\bn}$ connecting $S_0$ with $T$, and an even number of occurrences of $j$ in the walk $W_{\bm}W_{\bn}W_{\tilde\tau}$ connecting $S_0$ with $S$. It follows that $j\in J_T\setminus J_S$. Thus, $J_T=\gg_j(J_S)=J_S+\{j\}$, so $\bb$ is well-defined. Moreover, the above arguments show that
$$
P\tau=Q\quad\eq\quad \aa(P)\bb(\tau)=\aa(Q),
$$
for any $P,Q\in\SSS$ and $\tau\in\TTT$. It is clear that $\aa$ and $\bb$ are bijections, so $(\aa,\bb)$ is an isomorphism from $(\SSS,\TTT)$ onto $(\FFF,\TTT_G)$.
\end{proof}


Since cubical systems $(\SSS,\TTT)$ and $(\FFF,\TTT_G)$ of Theorem~\ref{RT} are isomorphic, their graphs are isomorphic to the graph $G$. The next result is obvious.

\begin{theorem} \label{quasi=cubical}
The graph of a cubical system is cubical. Conversely, any cubical graph $G$ defines a cubical system (a $G$-system).
\end{theorem}

\section{Stochastic token cubical systems} \label{S:stochastic}

Following~\cite{jF97} we consider a discrete stochastic process arising when random events result in occurrences of tokens in a finite cubical system $(\SSS,\TTT)$.

\begin{definition}
{\rm A quadruple $(\SSS,\TTT,\xi,\theta)$ is a {\em probabilistic token cubical system} if the following three conditions hold:
\roster
	\item[(i)] $(\SSS,\TTT)$ is a cubical system.
	\item[(ii)] $\xi:S\mapsto\xi(S)$ is a probability distribution (the {\em initial distribution}) on $\SSS$.
	\item[(iii)] $\theta:\tau\mapsto\theta_\tau$ is a probability distribution on $\TTT$ with $\theta_\tau>0$ for all tokens $\tau$ in $\TTT$.
\endroster
}
\end{definition}

Selecting an initial state according to the distribution $\xi$, and applying occuring tokens first to the initial state and then to its images under successive tokens, we obtain a Markov chain which we denote by $(\BS_n)$ where $n$ is the number of trials. The transition matrix $\BP$ of this chain is given by the equations
$$
p(S,V)=\begin{cases}
	\theta_\tau &\text{if $S\tau=V$,}\\
	0 &\text{otherwise,}
\end{cases}\quad\text{for $V\neq S$,}
$$
and
$$
p(S,S)=1-\sum_{V\in\SSS\setminus\{S\}}p(S,V).
$$
Note that $0<p(S,S)<1$, since, by Axiom [C2], for any state $S$ of the cubical system there is a token $\tau$ which is effective for $S$ with $\theta_\tau>0$.

The $n$-step transition probabilities are
\begin{align*}
p^{(n)}(S,V)&=\text{Pr}(\BS_n=V\mid \BS_0=S)\\
&=\sum_{(S_i)}\text{Pr}(\BS_n=S_n\mid \BS_{n-1}=S_{n-1},\ldots,\BS_0=S_0)\\
&=\sum_{(S_i)}p(S_0,S_1)p(S_1,S_2)\cdots p(S_{n-1},S_n),
\end{align*}
where the sums are taken over all $n$-tuples of states $(S_i)=(S_0,S_1,\ldots,S_n)$ with $S_0=S$ and $S_n=V$. Thus numbers $p^{(n)}(S,V)$ are entries of the matrix $\BP^n$.

\begin{lemma} \label{N-sequence}
Let $n\ge |\SSS|-1$. For any two states $S,V\in\SSS$ there is a sequence of states $S_0=S,S_1,\ldots,S_n=V$ such that, for any $0\le i<n$, the consecutive states $S_i$ and $S_{i+1}$ are either adjacent or equal. 
\end{lemma}

\begin{proof}
If $V=S$, we take $S_i=S$ for all $i$. Otherwise, by Axiom [C2], there is a stepwise effective message $\bm=\tau_1\ldots\tau_m$ producing $V$ from $S$. We may assume that the states $S_i$ produced by this message from $S$ are all distinct (take the shortest path in the graph of the cubical system). Then $m\le n$ and the $n$-tuple $(S_0,\ldots,S_m,S_m,\ldots,S_m)$ satisfies  conditions of the lemma.
\end{proof}

Let $S$ and $V$ be two states of the cubical system and $(S_i)$ be an $n$-tuple satisfying conditions of Lemma~\ref{N-sequence}. Then
$$
p(S_0,S_1)p(S_1,S_2)\cdots p(S_{n-1},S_n)>0.
$$
It follows that $p^{(n)}(S,V)>0$ for all $n\ge|\SSS|-1$ and $S,V\in\SSS$. Thus $(\BS_n)$ is a regular Markov chain.

Consider quantities $t(S)=\prod_{\tau\in\widehat S}\theta_\tau>0$. If $S$ and $V$ are two adjacent states, then, by Lemma~\ref{hatV-hatS}, $\widehat V\setminus\widehat S=\{\tau\}$ and $\widehat S\setminus\widehat V=\{\tilde\tau\}$ for some token $\tau$. Therefore,
$$
\frac{t(V)}{t(S)}=\frac{\prod_{\mu\in\widehat V}\theta_\mu}{\prod_{\mu\in\widehat S}\theta_\mu}=\frac{\theta_\tau}{\theta_{\tilde\tau}}=\frac{p(S,V)}{p(V,S)},
$$
so
$$
t(V)p(V,S)=p(S,V)t(S).
$$
Clearly, the last identity holds for all pairs of states $\{S,V\}$. Defining the probability distribution $\pi$ on $S$ by
$$
\pi(S)=\frac{t(S)}{\sum_{R\in\SSS}t(R)},
$$
we obtain
$$
\pi(V)p(V,S)=p(S,V)\pi(S).
$$
It follows that $(\BS_n)$ is a reversible regular Markov chain. Therefore, $\pi$ is its unique stationary distribution~\cite{jN97}. We established the following result (cf. Theorem 5.2 in ~\cite{jF97}):

\begin{theorem} \label{Markov}
The stochastic process $(\BS_n)$ is a regular Markov chain on the set of states $\SSS$. The unique asymptotic probability distribution $\pi$ on $\SSS$ is specified by
$$
\pi(S)=\frac{\prod_{\tau\in\widehat S}\theta_\tau}{\sum_{R\in\SSS}\prod_{\tau\in\widehat R}\theta_\tau}.
$$
\end{theorem}

\section{Examples} \label{S:examples}

We begin by introducing a class of finite $G$-systems that serves as a source of our examples (cf.~\cite{cD01}).

\begin{definition}
{\rm Let $\FFF$ be a family of subsets of a finite set $X$ with $|\FFF|\ge 2$. A set $S\in\FFF$ is said to be {\em downgradable} if there exists $x\in S$ such that $S\setminus\{x\}\in\FFF$. The family $\FFF$ itself is {\em downgradable} if all its nonminimal sets are downgradable.
Likewise, a set $S\in\FFF$ is said to be {\em upgradable} if there exists $x\in X\setminus S$ such that $S\cup\{x\}\in\FFF$. The family $\FFF$ itself is {\em upgradable} if all but its maximal sets are upgradable.
}
\end{definition}

It is clear that any downgradable family of sets containing the empty set is connected. Likewise, any upgradable family of subsets of $X$ containing the set $X$ itself is connected. Let $\FFF$ be any of such families. Then the induced subgraph $\langle\FFF\rangle$ of the cube $\HHH(X)$ is connected and therefore defines a cubical system (an $\langle\FFF\rangle$-system). 

\begin{example}
{\rm {\sl Comparability Graphs.} A simple finite graph $G=(X,E)$ is called a {\em comparability graph}~\cite{pF85} if there exists a partial order $P$ on $X$ such that
\begeq \label{comp graph}
\{x,y\}\in E \quad\eq\quad (x,y)\in P\quad\text{or}\quad (y,x)\in P.
\edeq
We denote $\mathfrak G$ the family of all comparability graphs on a fixed set $X$ and identify this family with the family of all sets of edges of comparability graphs on $X$. Note that $\mathfrak G$ contains the empty graph on $X$. It is known (see, for instance,~\cite{jD97}) that the family $\PPP$ of all partial orders on $X$ is well-graded and therefore is downgradable since it contains the empty partial order. As it can be easily seen this fact implies that the family $\mathfrak G$ is downgradable and therefore defines a cubical system.

Note that the wellgradedness property of the family $\PPP$ does not imply that $\mathfrak G$ is well-graded (see the graphs in Figure~\ref{comparability}).

{\begin{figure}[h!]
\centerline{\includegraphics{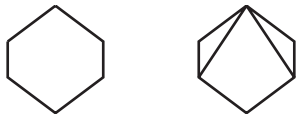}}
\caption{Two comparability graphs with 6 and 8 edges, respectively. The distance between the two edge sets is 2. There is no comparability graph on distance 1 from each of these two graphs.} \label{comparability} 
\end{figure}
}
}
\end{example}

\begin{example}
{\rm {\sl Interval and Indifference Graphs.} Interval and indifference graphs~\cite{pF85} are complements of comparability graphs arising from interval orders and semiorders, respectively, via relation~(\ref{comp graph}). As the families of all interval orders and all semiorders are well-graded~\cite{jD97} and both contain the empty relation, the respective families of interval and indifference graphs are upgradable and both contain the complete graph on $X$. Thus we can cast each of these two families as a cubical system.

Note that the same result can be obtained for any family of indifference graphs associated with partial orders satisfying so-called ``distinguishing property''~\cite{sO05}.
}
\end{example}

\begin{example}
{\rm {\sl Almost Connected Orders.} An {\em ac-order} (almost connected order)~\cite{cD01} is an asymmetric binary relation $R$ on a set $X$ satisfying the following condition:
$$
(x,y)\in R\;\;\text{and}\;\;(y,z)\in R\quad\imp\quad(x,w)\in R\;\;\text{or}\;\;(w,z)\in R
$$
for all $x$, $y$, $z$, $w$ in $X$. It is shown in~\cite{cD01} that the family $\AAA$ of all ac-orders on $X$ is both downgradable and upgradable and connected. We conclude that the family $\AAA$ can be cast as a cubical system. Note that Theorem~29 in~\cite{cD01} asserts that $\AAA$ is not well-graded if $|X|>4$.
}
\end{example}

We conclude this section with a simple example of an infinite cubical system.

\begin{example}
{\rm Let $\ZZZ^n$ be the graph of the $n$-dimensional integer lattice $\Zee^n$. It is not difficult to show that $\ZZZ^n$ is isometrically embeddable into some (infinite) cube $\HHH(X)$. Thus any connected subgraph $G$ of the graph $\ZZZ^n$ is cubical and therefore defines a $G$-system.
}
\end{example}

\section{Conclusion}

We have investigated algebraic and stochastic properties of cubical systems and shown that main results of media theory hold for cubical systems. The relations between families of media and cubical systems on a given set of states are indicated in the diagram shown below:

{\begin{figure}[h!]
\centerline{\includegraphics{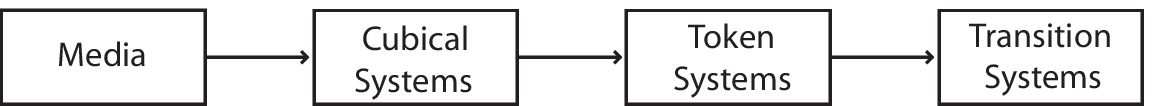}}
\end{figure}
}

The structural properties of message and state contents (Theorem~\ref{hatV-hatS}), together with the representation theorem (Theorem~\ref{RT}), reveal the binary nature of states in both media and cubical systems theories, which is also demonstrated by the `cubical' structure of the corresponding graphs (Theorem~\ref{quasi=cubical}). This characterization of states is crucial for the stochastic token theory (Theorem~\ref{Markov}). Because any subgraph of a cube is a disjoint union of connected cubical graphs, it is appropriate to say that cubical systems represent the most general case of token systems enjoying the binary structure of their states.

Our treatment of cubical systems as token systems rather than transition systems is motivated by examples in Section~\ref{S:examples} and connections with media theory. On the other hand, general methods of ``concurrency'' theory~\cite{gW94}, and especially ``geometric'' models for concurrency~\cite{vP91,eG95} could bring new elements to cubical token systems theory. In particular, a topological cubical complex can be associated with a cubical system in a natural way. Such complexes were used in the treatment of weak order families as media in~\cite{sO04,sO05}.

\end{document}